
\input amstex
\documentstyle{amsppt}
\magnification=1200
\catcode`\@=11
\redefine\logo@{}
\catcode`\@=13
\pageheight{19cm}

\define \bn{\Bbb N}
\define \bz{\Bbb Z}
\define \bq{\Bbb Q}
\define \br{\Bbb R}
\define \bc{\Bbb C}

\define\Da{{\Cal D}}
\define \M{{\Cal M}}

\define\Ka{{\Cal K}}
\define \E{{\Cal E}}

\define\rk{\text{rk}~}


\define\pr{\text{pr}}

\define\Exc{\text{Exc}}
\define\0o{{\overline 0}}
\define\1o{{\overline 1}}

\TagsOnRight


\topmatter

\title
On a classical correspondence between K3 surfaces II
\endtitle

\author
Carlo Madonna and
Viacheslav V. Nikulin \footnote{Supported by Russian Fund of
Fundamental Research (grant N 00-01-00170).\hfill\hfill}
\endauthor

\address
Dipartimento di Matematica, Universit\`a degli studi di Ferrara, Italia
\endaddress
\email
madonna\@mat.uniroma3.it
\endemail

\address
Deptm. of Pure Mathem. The University of Liverpool, Liverpool
L69 3BX, UK;
\vskip1pt
Steklov Mathematical Institute,
ul. Gubkina 8, Moscow 117966, GSP-1, Russia
\endaddress
\email
vnikulin\@liv.ac.uk\ \
slava\@nikulin.mian.su
\endemail

\dedicatory
To memory of Andrei Nikolaevich Tyurin 
\enddedicatory

\abstract
Let $X$ be a K3 surface and $H$ a primitive polarization of degree 
$H^2=2a^2$, $a>1$. The moduli space of sheaves over $X$ 
with the isotropic Mukai vector $(a,H,a)$ is again a K3 surface $Y$ 
which is endowed by a natural $nef$ element $h$ with $h^2=2$. We give 
necessary and sufficient conditions in terms of Picard lattices $N(X)$ 
and $N(Y)$ when $Y\cong X$, generalising our results in \cite{4} for $a=2$.

In particular, we show that $Y\cong X$ if for at least one 
$\alpha =\pm 1$ there exists 
$h_1\in N(X)$ such that $h_1^2= 2\alpha a$,    
$H\cdot h_1\equiv 0\mod a$, and the primitive sublattice 
$[H,h_1]_{pr} \subset N(X)$ contains $x$ such that 
$x\cdot H=1$.  

We also show that all divisorial conditions on moduli 
of $(X,H)$ (i.e. for Picard number 2) which imply $Y\cong X$ 
and $H\cdot N(X)=\bz$ are labelled by pairs 
$(d,\,\pm \mu)$ where $d\in \bn$, $\pm \mu \subset (\bz/2a^2)^\ast$ 
such that $d\equiv \mu^2\mod 4a^2$ and at least for one of $\alpha =\pm 1$  
the equation $p^2-dq^2=4a/\alpha$ has an integral solution $(p,q)$ 
with $p\equiv \mu q\mod 2a$. For each 
such $\pm \mu$ and $\alpha$, the number of $d$ and the corresponding 
divisorial conditions is infinite. Some of these conditions were found 
(in different form) by A.N. Tyurin in 1987. 
\endabstract

\rightheadtext
{On correspondences between K3 surfaces II}
\leftheadtext{C. Madonna and V.V. Nikulin}
\endtopmatter

\document

\head
0. Introduction
\endhead

Let $X$ be a K3 surface with a primitive
polarization $H$ of degree $H^2=2a^2$. Let $Y$ be the
moduli of sheaves over $X$ with the isotropic Mukai vector $v=(a,H,a)$
(see \cite{6}, \cite{7}). The $Y$ is a K3 surface which is endowed by 
a natural $nef$ element $h$ with $h^2=2$. It is isogenous to $X$ in 
the sense of Mukai. 

\proclaim{Question 1} When is $Y$ isomorphic to $X$?
\endproclaim 

We want to answer this question in terms of Picard lattices
$N(X)$ and $N(Y)$ of
$X$ and $Y$. Then our question reads as follows:

\proclaim{Question 2} Assume that $N$ is a hyperbolic lattice,
$\widetilde{H}\in N$ a primitive element with square $2a^2$. What are
conditions on $N$ and $\widetilde{H}$ such that for any K3 surface
$X$ with Picard lattice $N(X)$ and a primitive polarization $H\in N(X)$
of degree $2a^2$ the corresponding K3 surface $Y$ is isomorphic to $X$
if the the pairs of lattices $(N(X), H)$ and
$(N, \widetilde{H})$ are isomorphic as abstract lattices with
fixed elements?

In other words, what are conditions on $(N(X), H)$ as an abstract lattice
with a primitive vector $H$ with $H^2=2a^2$ which are sufficient for
$Y$ to be isomorphic to $X$ and are necessary
if $X$ is a general K3 surface with Picard lattice $N(X)$?
\endproclaim 

In \cite{4} we answered this question when $a=2$. 
Here we give an answer for any $a \geq 2$. For odd $a$, we additionally 
assume that $H\cdot N(X)=\bz$. For even $a$, that is valid, 
if $Y\cong X$. The answer is given in Theorem 2.2.2 and
also Proposition 2.2.1. 

In particular, if the Picard number
$\rho (X)=\rk N(X) \ge 12$, the result is very simple:
$Y\cong X$, if and only if either there exists
$x\in N(X)$ such that $x\cdot H=1$ or $a$ is odd and 
there exists $x\in N(X)$ such that $x\cdot H=2$. 
This follows from results of Mukai \cite{7} and also  \cite{9}, \cite{10}.

The polarized K3 surfaces $(X,H)$ with
$\rho (X)=2$ are especially interesting. 
Really, it is well-known that the moduli space of 
polarized K3 surfaces of any even degree $H^2$ is $19$-dimensional. 
If $X$ is general, i. e. $\rho(X)=1$, then the surface
$Y$ cannot be isomorphic to $X$
because $N(X)=\bz H$ where $H^2=2a^2$, $a>1$, and $N(X)$ does not 
have elements with square $2$ which is necessary if $Y\cong X$. 
Thus, if $Y\cong X$, then $\rho (X)\ge 2$, and $X$ belongs 
to a codimension $\rho (X)-1$ submoduli space of K3 surfaces which is a
{\it divisor in the 19-dimensional moduli space $\M$  
(up to codimension 2).} To describe connected components of the divisor, 
it is equivalent to describe Picard lattices $N(X)$ of the surfaces
$X$ with fixed $H\in N(X)$ such that $\rho (X)=\rk N(X)=2$ 
and a general K3 surfaces $X$ with Picard lattice $N(X)$ has 
$Y\cong X$. 

The pair $(N(X),H)$ with $\rho =\rk N(X)=2$ and $H\cdot N(X)=\bz$ 
is defined up to isomorphisms by $d=-\det N(X)$  $>0$ 
(it defines the Picard lattice $N(X)$ up to isomorphisms, if $Y\cong X$) 
and by the invariant 
$\pm \mu =\{\mu,-\mu\}\subset (\bz/2a^2)^\ast$ (this is the invariant of 
the primitive vector $H\in N(X)$) such that 
$\mu^2 \equiv d\mod 4a^2$. See Proposition 3.1.1 about 
definition of $\mu$. 
We show that for a general $X$ with such $N(X)$ 
we have $Y\cong X$ and for odd $a$ additionally $H\cdot N(X)=\bz$, 
if and only if at least for one $\alpha=\pm 1$ 
there exists integral $(p,q)$ such that 
$$
p^2-dq^2=4a/\alpha\ \text{and}\ p\equiv \mu q\mod 2a .
\tag{0.1}
$$
For each such $\pm \mu$ and $\alpha$ the set 
$\Da_{\alpha}^\mu$ of $d$ having such solution $(p,q)$ and 
$d\equiv \mu^2\mod {4a^2}$ is infinite since it 
contains the infinite subset 
$$
\{(\mu+2ta/\alpha)^2-4a/\alpha>0\ |\ t\mu \equiv 1\mod a \}
\tag{0.2}
$$
(put $q=1$ to \thetag{0.1}). Thus, the set of possible 
divisorial conditions on moduli of $(X,H)$ 
which imply $Y\cong X$ and $H\cdot N(X)=\bz$ is 
labelled by the set of pairs $(d, \pm \mu)$ described above, 
and it is infinite. 

Some of infinite series of divisorial 
conditions on moduli of $X$ which imply $Y\cong X$ were found 
by A.N. Tyurin in \cite{17} --- \cite{19}. He found, in different 
form, infinite series \thetag{0.2} for $\alpha=-1$ and any $\pm \mu$  
(if not other ones).  

Surprisingly, solutions $(p,q)$ of \thetag{0.1} can be 
interpreted as elements of Picard lattices of $X$ and $Y$. 
We get the following simple sufficient condition on $(X,H)$ 
which implies $Y\cong X$. 
It seems, many known examples when it happened that $Y\cong X$ 
(e. g. see \cite{2}, \cite{8}, \cite{17}) 
follow from this condition. This is one of the
main results of the paper, and we want to formulate it exactly      
(a similar statement can be also formulated in terms of $Y$). 

\proclaim{Theorem}
Let $X$ be a K3 surface with a primitive polarization $H$ of degree $2a^2$, 
$a \geq 2$. Let $Y$ be the moduli space of sheaves on $X$ with the 
Mukai vector $v=(a,H,a)$.

Then $Y\cong X$, if at least for one $\alpha=\pm 1$  
there exists $h_1\in N(X)$ such that
$$
(h_1)^2= 2\alpha a,\ \  h_1\cdot H\equiv 0 \mod a,
$$
and the primitive sublattice  $[H,h_1]_{\pr}\subset N(X)$ generated by 
$H$, $h_1$, contains $x$ such that $x\cdot H=1$. 

These conditions are necessary to have $Y\cong X$ and 
$H\cdot N(X)=\bz$ for $a$ odd, if either $\rho(X)=1$, 
or $\rho(X)=2$ and $X$ is a general K3 surface with its Picard lattice.
\endproclaim

From our point of view, this statement is also very interesting
because some elements $h_1$ of the Picard lattice $N(X)$ with negative
square $(h_1)^2$ get a very clear geometrical meaning (when $\alpha<0$). 
For K3 surfaces this is well-known only for elements 
$\delta$ of the Picard lattice $N(X)$ with negative square 
$\delta^2=-2$: then $\delta$ or $-\delta$ is effective.

\smallpagebreak 

As for the case $a=2$,
the fundamental tool to get the results above is the Global
Torelli Theorem for K3 surfaces \cite{12} and results of Mukai 
\cite{6}, \cite{7}. Using the results of Mukai, we can calculate periods of 
$Y$ using periods of $X$; by the Global Torelli Theorem \cite{12}, 
we can find out if $Y$ is isomorphic to $X$.

\head
1. Preliminary notations and results about lattices and K3 surfaces
\endhead

\subhead
1.1. Some notations about lattices
\endsubhead
We use notations and terminology from \cite{10} about lattices,
their discriminant groups and forms. A {\it lattice} $L$ is a
non-degenerate integral symmetric bilinear form. I. e. $L$ is a free
$\bz$-module equipped with a symmetric pairing $x\cdot y\in \bz$ for
$x,\,y\in L$, and this pairing should be non-degenerate. We denote
$x^2=x\cdot x$. The {\it signature} of $L$ is the signature of the
corresponding real form $L\otimes \br$.
The lattice $L$ is called {\it even}
if $x^2$ is even for any $x\in L$. Otherwise, $L$ is called {\it odd}.
The {\it determinant} of $L$ is defined to be $\det L=\det(e_i\cdot e_j)$
where $\{e_i\}$ is some basis of $L$. The lattice $L$ is {\it unimodular}
if $\det L=\pm 1$.
The {\it dual lattice} of $L$ is
$L^\ast=Hom(L,\,\bz)\subset L\otimes \bq$. The
{\it discriminant group} of $L$ is $A_L=L^\ast/L$. It has the order
$|\det L|$. The group $A_L$ is equipped with the
{\it discriminant bilinear form} $b_L:A_L\times A_L\to \bq/\bz$
and the {\it discriminant quadratic form} $q_L:A_L\to \bq/2\bz$
if $L$ is even. To get this forms, one should extend the form of $L$ to
the form on the dual lattice $L^\ast$ with values in $\bq$.

For $x\in L$, we shall consider
the invariant $\gamma(x)\ge 0$ where
$$
x\cdot L=\gamma (x)\bz .
\tag{1.1.1}
$$
Clearly, $\gamma (x)|x^2$ if $x\not=0$.

We denote by $L(k)$ the lattice obtained from a lattice $L$ by
multiplication of the form of $L$ by $k\in \bq$.
The orthogonal sum of lattices $L_1$ and $L_2$ is denoted by
$L_1\oplus L_2$.
For a symmetric integral matrix
$A$, we denote by $\langle A \rangle$ a lattice which is given by
the matrix $A$ in some bases. We denote
$$
U=\left(
\matrix
0&1\\
1&0
\endmatrix
\right).
\tag{1.1.2}
$$
Any even unimodular lattice of the signature $(1,1)$ is isomorphic to
$U$.

An embedding $L_1\subset L_2$ of lattices is called {\it primitive}
if $L_2/L_1$ has no torsion.
We denote by $O(L)$, $O(b_L)$ and $O(q_L)$ the automorphism groups of
the corresponding forms.
Any $\delta\in L$ with $\delta^2=-2$ defines
a reflection $s_\delta\in O(L)$ which is given by the formula
$$
x\to x+(x\cdot \delta)\delta,
$$
$x\in L$. All such reflections generate
the {\it 2-reflection group} $W^{(-2)}(L)\subset O(L)$.

\subhead
1.2. Some notations about K3 surfaces
\endsubhead
Here we remind some basic notions and results about K3 surfaces,
e. g. see \cite{12}, \cite{13}, \cite{14}.
A K3 surface $S$ is a non-singular projective algebraic surface over
$\bc$ such that its canonical class $K_S$ is zero and the irregularity
$q_S=0$. We denote by $N(S)$ the {\it Picard lattice} of $S$ which is
a hyperbolic lattice with the intersection pairing
$x\cdot y$ for $x,\,y\in N(S)$. Since the canonical class $K_S=0$,
the space $H^{2,0}(S)$ of 2-dimensional holomorphic differential
forms on $S$ has dimension one over $\bc$, and
$$
N(S)=\{x\in H^2(S,\bz)\ |\ x\cdot H^{2,0}(S)=0\}
\tag{1.2.1}
$$
where $H^2(S,\bz)$ with the intersection pairing is a
22-dimensional even unimodular lattice of signature
$(3,19)$. The orthogonal lattice $T(S)$ to $N(S)$ in $H^2(S,\bz)$ is called
the {\it transcendental lattice of $S$.} We have
$H^{2,0}(S)\subset T(S)\otimes \bc$. The pair $(T(S), H^{2,0}(S))$ is
called the {\it transcendental periods of $S$}.
The {\it Picard number} of $S$ is
$\rho(S)=\rk N(S)$. A non-zero element $x\in N(S)\otimes \br$ is
called {\it nef} if $x\not=0$ and $x\cdot C\ge 0$ for any effective
curve $C\subset S$. It is known that an element $x\in N(S)$ is ample
if $x^2>0$, $x$ is $nef$, and the orthogonal complement
$x^\perp$ to $x$ in  $N(S)$ has no elements with square $-2$.
For any element $x\in N(S)$ with $x^2\ge 0$, there exists a reflection
$w\in W^{(-2)}(N(S))$ such that the element
$\pm w(x)$ is nef; it then is ample, if $x^2>0$ and $x^\perp$ had no
elements with square $-2$ in $N(S)$.

We denote by $V^+(S)$ the light cone of $S$, which is the half-cone of
$$
V(S)=\{x\in N(S)\otimes \br\ |\ x^2>0\ \}
\tag{1.2.2}
$$
containing a polarization of $S$. In particular, all $nef$ elements
$x$ of $S$ belong to $\overline{V^+(S)}$:
one has $x\cdot V^+(S)>0$ for them.

The reflection group $W^{(-2)}(N(S))$ acts in $V^+(S)$ discretely,
and its  fundamental
chamber is the closure $\overline{\Ka(S)}$ of the K\"ahler cone
$\Ka(S)$ of $S$. It is the same as the set of all $nef$ elements of $S$.
Its faces are orthogonal to the set $\Exc(S)$ of all exceptional curves $r$
on $S$ which are non-singular rational curves $r$ on $S$ with $r^2=-2$.
Thus, we have
$$
\overline{\Ka (S)}=\{0\not=x\in \overline{V^+(S)}\ |
\ x\cdot \Exc(S)\ge 0\,\}.
\tag{1.2.3}
$$

\head
2. General results on the Mukai correspondence between K3 surfaces
with primitive polarizations of degrees $2a^2$ and 2 which gives 
isomorphic K3's
\endhead

\subhead
2.1. The correspondence
\endsubhead
Let $X$ be
a K3 surface with a primitive
polarization $H$ of degree $2a^2$, $a>0$. 
Let $Y$ be the moduli space of (coherent) sheaves $\E$ on $X$ with the 
isotropic Mukai vector $v=(a,H,a)$. It means that $\rk \E=a$, 
$\chi (\E)=2a$ and $c_1(\E)=H$. 
Let
$$
H^\ast(X,\bz)=H^0(X,\bz)\oplus H^2(X,\bz)\oplus H^4(X,\bz)
\tag{2.1.1}
$$
be the full cohomology lattice of $X$ with the Mukai product
$$
(u,v)=-(u_0\cdot v_2+u_2\cdot v_0)+u_1\cdot v_1
\tag{2.1.2}
$$
for $u_0,v_0\in H^0(X,\bz)$, $u_1,v_1\in H^2(X,\bz)$,
$u_2,v_2\in H^4(X,\bz)$. We naturally identify
$H^0(X,\bz)$ and $H^4(X,\bz)$ with $\bz$. Then the Mukai product is
$$
(u,v)=-(u_0v_2+u_2v_0)+u_1\cdot v_1.
\tag{2.1.3}
$$
The element
$$
v=(a,H,a)=(a,H,\chi-a)\in H^\ast(X,\bz)
\tag{2.1.4}
$$
is isotropic, i.e. $v^2=0$. In this case,  
Mukai showed \cite{6}, \cite{7} that $Y$ is a K3 surface, and 
one has the natural identification
$$
H^2(Y,\,\bz)= (v^\perp/\bz v)
\tag{2.1.5}
$$
which also gives the isomorphism of the Hodge structures of $X$ and $Y$.
The element $h=(-1,0,1)\in v^\perp$ has square $h^2=2$, $h\mod \bz v$ belongs
to the Picard lattice $N(Y)$ of $Y$ and is $nef$. 
See \cite{5}, \cite{13} and \cite{15} about geometry of $(Y,h)$. 
For a general $X$, the K3 surface $Y$ is a double plane. 

We want to answer Question 2 which we exactly formulated in 
Introduction: Using $N(X)$, say when $Y\cong X$.

\subhead
2.2. Formulation of general results
\endsubhead
We use notations from Sect. 2.1. Thus, we assume that $X$ is a K3 
surface with a primitive polarization $H$ with $H^2=2a^2$ where $a>1$.  
The following statement follows from results of Mukai \cite{7} 
and some results from \cite{10}. It is standard and 
well-known.

\proclaim{Proposition 2.2.1} 
If $Y$ is isomorphic to $X$, then either $\gamma (H)=1$ or 
$a$ is odd and $\gamma (H)=2$ for $H\in N(X)$  
(see \thetag{1.1.1}).

Assume that either $\gamma (H)=1$ or $a$ is odd and $\gamma (H)=2$ for 
$H\in N(X)$.  
Then the Mukai identification \thetag{2.1.5}
canonically identifies the transcendental periods
$(T(X), H^{2,0}(X))$ and 
$(T(Y), H^{2,0}(Y))$. It follows that
the Picard lattices $N(Y)$ and $N(X)$ have the same genus.
In particular, $N(Y)$ is isomorphic to $N(X)$ if the genus of $N(X)$
contains only one class.  If the genus of $N(X)$ contains only one class,
then $Y$ is isomorphic to $X$, if additionally the canonical homomorphism
$O(N(X))\to O(q_{N(X)})$ is epimorphic.
Both these conditions are valid (in particular, $Y\cong X$),
if $\rho(X)\ge 12$.
\endproclaim

From now on we assume that $\gamma(H)=1$ in $N(X)$, which  
is automatically valid for even $a$, if $Y\cong X$.

Calculations below are valid for an arbitrary K3 surface $X$ and
a primitive vector $H\in N(X)$ with $H^2=2a^2$, $a>0$ and $\gamma (H)=1$.
Let $K(H)=H^\perp_{N(X)}$ be the orthogonal complement to $H$ in $N(X)$.
Set $H^\ast=H/2a^2$. Then any element $x\in N(X)$ can be written as
$$
x=nH^\ast+k^\ast
\tag{2.2.1}
$$
where $n\in \bz$ and $k^\ast\in K(H)^\ast$, because
$$
\bz H\oplus K(H)\subset N(X)\subset
N(X)^\ast\subset \bz H^\ast\oplus K(H)^\ast. 
$$
Since $\gamma(H)=1$,
the map $nH^\ast+[H] \mapsto k^\ast +K(H)$ gives an isomorphism of
the groups
$\bz/2a^2 \cong [H^\ast]/[H]\cong [u^\ast+K(H)]/K(H)$ where
$u^\ast+K(H)$ has order $2a^2$ in $A_{K(H)}=K(H)^\ast/K(H)$.
It follows,
$$
N(X)=[\bz H, K(H), H^\ast+u^\ast].
\tag{2.2.2}
$$
The element $u^\ast$ is defined canonically mod $K(H)$. Since
$H^\ast+u^\ast$ belongs to the even lattice $N(X)$, it follows
$$
(H^\ast+u^\ast)^2={1\over {2a^2}}+{u^\ast}^2 \equiv 0 \mod 2.
\tag{2.2.3}
$$
Let $\overline{H^\ast}=H^\ast \mod [H]\in [H^\ast]/[H]\cong \bz/2a^2$
and $\overline{k^\ast}=k^\ast\mod K(H)\in A_{K(H)}=K(H)^\ast/K(H)$.
We then have
$$
N(X)/[H,K(H)]=(\bz/2a^2)(\overline{H^\ast} +\overline{u^\ast})\subset
(\bz/2a^2)\overline{H^\ast}+K(H)^\ast/K(H).
\tag{2.2.4}
$$
Also $N(X)^\ast\subset \bz H^\ast+K(H)^\ast$ since
$H+K(H)\subset N(X)$,
and for $n\in \bz$, $k^\ast \in K(H)^\ast$ we have
$x=nH^\ast+k^\ast \in N(X)^*$ if and only if
$$
(nH^\ast+k^\ast)\cdot (H^\ast+u^\ast)=
{n\over 2a^2}+k^\ast\cdot u^\ast\in \bz.
$$
It follows,
$$
N(X)^\ast=
\{nH^\ast+k^\ast\ |\ n\in \bz,\ k^\ast \in K(H)^\ast,\
n\equiv -2a^2 \ (k^\ast\cdot u^\ast) \mod 2a^2 \}\subset
$$
$$
\subset \bz H^\ast+K(H)^\ast ,
\tag{2.2.5}
$$
and
$$
\split
N(X)^\ast/[H,K(H)]&=
\{-2a^2(\overline{k^\ast}\cdot \overline{u^\ast})\, \overline{H^\ast}+
\overline{k^\ast}\}\ |\ \overline{k^\ast} \in A_{K(H)}\}
\subset \\
&\subset (\bz/2a^2) \overline{H^\ast}+A_{K(H)}.
\endsplit
\tag{2.2.6}
$$
We introduce {\it the characteristic map of the polarization $H$}
$$
\kappa(H):K(H)^\ast \to A_{K(H)}/(\bz/2a^2)(u^\ast+K(H))\to A_{N(X)}
\tag{2.2.7}
$$
where for $k^\ast \in K(H)^\ast$ we have
$$
\kappa(H)(k^\ast)=-2a^2(k^\ast\cdot u^\ast) H^\ast+k^\ast + N(X)\in
A_{N(X)}.
\tag{2.2.8}
$$
It is epimorphic, its kernel is $(\bz/2a^2)(u^\ast+K(H))$, and it gives
the canonical isomorphism
$$
\overline{\kappa(H)}:A_{K(H)}/(\bz/2a^2)(u^\ast+K(H))
\cong A_{N(X)}.
\tag{2.2.9}
$$
For the corresponding discriminant forms we have
$$
\kappa(k^\ast)^2 \mod 2=(k^\ast)^2+2a^2(k^\ast\cdot u^\ast )^2\mod 2.
\tag{2.2.10}
$$

Now we can formulate our main result:

\proclaim{Theorem 2.2.2} The surface $Y$ is isomorphic to $X$
if the following conditions (a), (b), (c) are valid:

(a) $\gamma(H)=1$ for $H\in N(X)$;

(b) there exists $\widetilde{h}\in N(X)$ with
$\widetilde{h}^2=2$, $\gamma(\widetilde{h})=1$ and such that
there exists an embedding
$$
f:K(H)\to K(\widetilde{h})
$$
of negative definite lattices such that
$$
K(\widetilde{h})=[f(K(H)), 2a f(u^\ast)], \
w^\ast+K(\widetilde{h})=a f(u^\ast)+K(\widetilde{h});
$$

(c) the dual to $f$ embedding
$f^\ast:K(\widetilde{h})^\ast \to K(H)^\ast$ commutes
(up to multiplication by $\pm 1$) with the characteristic maps
$\kappa(H)$ and $\kappa(\widetilde{h})$, i. e.
$$
\kappa(\widetilde{h})(k^\ast)=\pm \kappa(H)(f^\ast(k^\ast))
\tag{2.2.11}
$$
for any $k^\ast \in K(\widetilde{h})^\ast$.

The conditions (a), (b) and (c) are necessary if for odd $a$ additionally 
$\gamma (H)=1$ for $H\in N(X)$, and  $\rk N(X)\le 19$,  
and $X$ is a general K3 surface with 
the Picard lattice $N(X)$ in the following sense: 
the automorphism group of the transcendental periods $(T(X),H^{2,0}(X))$ 
is $\pm 1$. (Remind that $Y\cong X$ if $\rk N(X)=20$.)
\endproclaim

\subhead
2.3. Proofs
\endsubhead
Let us denote by $e_1$ the canonical generator of $H^0(X,\bz)$ and
by $e_2$ the canonical generator of $H^4(X,\bz)$. They generate
the sublattice $U$ in $H^\ast (X,\bz)$.
Consider the Mukai vector $v=ae_1+H+ae_2=(a,H,a)$. We have
$$
N(Y)=v\,^\perp _{U\oplus N(X)}/\bz v.
\tag{2.3.1}
$$
Let us calculate $N(Y)$. Let $K(H)=H^\perp_{N(X)}$. Then we have
embedding of lattices of finite index
$$
\bz H\oplus K(H)\subset N(X)\subset N(X)^\ast\subset
\bz H^\ast \oplus K(H)^\ast
\tag{2.3.2}
$$
where $H^\ast =H/2a^2$.
We have the orthogonal decomposition up to finite index
$$
U\oplus \bz H\oplus K(H)\subset U\oplus N(X)\subset
U\oplus \bz H^\ast \oplus K(H)^\ast .
\tag{2.3.3}
$$
Let $s=x_1e_1+x_2e_2+yH^\ast+z^\ast\in v^\perp_{U\oplus N(X)}$,
$z^\ast \in K(H)^\ast$.
Then $-ax_1-ax_2+y=0$ since $s\in v^\perp$ and hence $(s,v)=0$.
Thus, $y=ax_1+ax_2$ and
$$
s=x_1e_1+x_2e_2+a(x_1+x_2)H^\ast + z^\ast .
\tag{2.3.4}
$$
Here $s\in U\oplus N(X)$ if and only if $x_1, x_2 \in \bz$ and
$a(x_1+x_2)H^\ast+z^\ast\in N(X)$. This orthogonal complement
contains
$$
[\bz v, K(H), \bz h]
\tag{2.3.5}
$$
where $h=-e_1+e_2$, and this is a sublattice of finite index
in $(v^\perp)_{U\oplus N(X)}$. The generators $v$, generators of $K(H)$
and $h$ are free, and we can rewrite $s$ above using these generators
with rational coefficients as follows:
$$
s={-x_1+x_2\over 2}h+{x_1+x_2\over 2a}v+z^\ast,
\tag{2.3.6}
$$
where $a(x_1+x_2)H^\ast+z^\ast \in N(X)$.
Equivalently, for $h^\ast=h/2$, 
$$
s=x_1'h^\ast+x_2'{v\over 2a}+z^\ast,
\tag{2.3.7}
$$
where $x_1',x_2'\in \bz$, $z^\ast \in K(H)^\ast$, $x_1'\equiv x_2'\mod 2$,
and $ax_2'H^\ast+z^\ast\in N(X)$.

From these calculations, we get 

\proclaim{Claim}
Assume that $\gamma (H)=1$. Then
$$
N(X)=[H,K(H),{H\over 2a^2}+u^\ast],
\tag{2.3.8}
$$
$$
N(Y)=[h, K(h)=[K(H),2a u^\ast],{h\over 2}+a u^\ast]=
[h, K(h),{h\over 2}+w^\ast],
\tag{2.3.9}
$$
where $u^\ast+K(H)$ has order $2a^2$ in $A_{K(H)}$, $w^\ast=a u^\ast$,
$K(h)=[K(H),2w^\ast=2a u^\ast]$. Here we agreed notations with Sect. 2.2.
We have $\det N(X)=\det K(H)/2a^2$ and $\det N(Y)=\det K(h)/2$ (in particular, 
$\gamma (h)=1$). Thus, $\det N(X)=\det N(Y)$ for this case, since 
$\det K(h)=\det K(H)/a^2$. 
We can formally put here $h={H\over a}$ since
$h^2=\left({H\over a}\right)^2=2$.
\endproclaim

From the claim, we get 

\proclaim{Lemma 2.3.1} For Mukai identification \thetag{2.1.5},
the sublattice $T(X)\subset T(Y)$
has index $1$, if $\gamma(H)=1$ for $H\in N(X)$.
\endproclaim

\demo{Proof} Really, since $H\in N(X)$, $T(X)\perp N(X)$ and
$T(X)\cap \bz v=\{0\}$, the Mukai identification \thetag{2.1.5} gives
an embedding
$T(X)\subset T(Y)$. We then have $\det T(Y)=\det T(X)/[T(Y):T(X)]^2$.
Moreover, $|\det T(X)|=|\det N(X)|$ and $|\det T(Y)|=|\det N(Y)|$
because the transcendental and the Picard lattice are orthogonal
complements to each other in a unimodular lattice $H^2(\ast,\bz)$.
By \thetag{2.3.8} and \thetag{2.3.9} we get the statement.
\enddemo

The statement of Lemma 2.3.1 is a particular
case of the general statement by Mukai \cite{7} that 
$$
[T(Y):T(X)]=q=\min |v\cdot x|
\tag{2.3.10}
$$
for all $x\in H^0(X,\bz)\oplus N(X)\oplus H^4(X,\bz)$ such that
$v\cdot x\not=0$. For our Mukai vector $v=(a,\,H,\,a)$ we obviously 
get that $q=1$, if and only if for $H\in N(X)$ either 
$\gamma (H)=1$ or $\gamma (H)=2$ and $a$ is odd. 
Thus, for $a$ even we have the only case: $\gamma (H)=1$.

\demo{Proof of Proposition 2.2.1} By \thetag{2.3.10}, if $Y\cong X$, 
we have either $\gamma (H)=1$ or $\gamma (H)=2$ and $a$ is odd. 

Assume that $\gamma (H)=1$ or $\gamma (H)=2$ and $a$ is odd. 
Then $T(X)=T(Y)$ for the Mukai identification \thetag{2.1.5}.
By the discriminant forms technique (see \cite{10}),
then the discriminant
quadratic forms $q_{N(X)}=-q_{T(X)}$ and $q_{N(Y)}=-q_{T(Y)}$ are isomorphic.
Thus, lattices $N(X)$ and $N(Y)$ have the same signatures and
discriminant quadratic forms. It follows (see \cite{10}) that they have
the same genus: $N(X)\otimes \bz_p\cong N(Y)\otimes \bz_p$ for any prime
$p$ and the ring of $p$-adic integers $\bz_p$.
Additionally, assume that either the genus of $N(X)$ or the
genus of $N(Y)$ contains only one class. Then $N(X)$ and
$N(Y)$ are isomorphic.

If additionally the canonical homomorphism $O(N(X))\to O(q_{N(X)})$
(equivalently, $O(N(Y))\to O(q_{N(Y)})$) is epimorphic, then
the Mukai identification $T(X)=T(Y)$ can be extended to give an
isomorphism $\phi:H^2(X,\bz)\to H^2(Y,\bz)$ of cohomology lattices.
The Mukai identification is identical on $H^{2,0}(X)=H^{2,0}(Y)$.
Multiplying $\phi$ by $\pm 1$ and by elements of the reflection
group $W^{(-2)}(N(X))$, if necessary, we
can assume that $\phi(H^{2,0}(X))=H^{2,0}(Y)$ and $\phi$ maps the
K\"ahler cone of $X$ to the K\"ahler cone of $Y$. By global Torelli
Theorem for K3 surfaces \cite{12}, $\phi$ is then defined 
by an isomorphism of K3 surfaces $X$ and $Y$.

If $\rho (X)\ge 12$, by \cite{10} Theorem 1.14.4, the primitive
embedding of $T(X)=T(Y)$ into the cohomology lattice $H^2(X,\bz)$ of
K3 surfaces is unique up to automorphisms of the lattice $H^2(X,\bz)$.
Like above, it then follows that $X$ is isomorphic to $Y$.
The given proof of Proposition 2.2.1 is standard and well-known. 
\enddemo

\demo{Proof of Theorem 2.2.2} 
Assume that $\gamma (H)=1$. The Mukai 
identification then gives the canonical identification
$$
T(X)=T(Y).
\tag{2.3.11}
$$
Thus, it gives the canonical identifications
$$
\split
A_{N(X)}=N(X)^\ast/N(X)=&(U\oplus N(X))^\ast/(U\oplus N(X))=
T(X)^\ast/T(X)=A_{T(X)}\\
=A_{T(Y)}=T(Y)^\ast/T(Y)=&N(Y)^\ast /N(Y)=A_{N(Y)}.
\endsplit
\tag{2.3.12}
$$
Here $A_{N(X)}=N(X)^\ast/N(X)=(U\oplus N(X))^\ast/(U\oplus N(X))$ because
$U$ is unimodular,
$(U\oplus N(X))^\ast/(U\oplus N(X))=T(X)^\ast/T(X)=A_{T(X)}$ because
$U\oplus N(X)$ and $T(X)$ are orthogonal complements to each other in
the unimodular lattice $H^\ast(X,\bz)$. Here \linebreak
$A_{T(Y)}=T(Y)^\ast/T(Y)=N(Y)^\ast/ N(Y)=A_{N(Y)}$
because $T(Y)$ and $N(Y)$ are
orthogonal complements to each other in the unimodular lattice
$H^2(Y,\bz)$.
E. g. the identification
$(U\oplus N(X))^\ast/(U\oplus N(X))=T(X)^\ast/T(X)=A_{T(X)}$ is
given by the canonical correspondence
$$
x^\ast+(U\oplus N(X))\to t^\ast +T(X)
\tag{2.3.13}
$$
if $x^\ast \in (U\oplus N(X))^\ast$, $t^\ast \in T(X)^\ast$ and
$x^\ast+t^\ast \in H^\ast(X,\bz)$.

By \thetag{2.3.9}, we also have the canonical embedding of lattices
$$
K(H)\subset K(h)=[K(H),2a u^\ast].
\tag{2.3.14}
$$

We have the key statement:

\proclaim{Lemma 2.3.2} Assume that $\gamma (H)=1$. The canonical embedding 
\thetag{2.3.14} (it is given by \thetag{2.3.9}) 
$K(H)\subset K(h)$ of lattices,
and the canonical identification
$A_{N(X)}=A_{N(Y)}$ (given by \thetag{2.3.12}) agree with the
characteristic homomorphisms $\kappa(H):K(H)^\ast\to A_{N(X)}$
and $\kappa(h):K(h)^\ast \to A_{N(Y)}$, i.e.
$\kappa(h)(k^\ast)=\kappa (H)(k^\ast)$ for any $\kappa^\ast \in
K(h)^\ast\subset K(H)^\ast$ (this embedding is dual to \thetag{2.3.14}).
\endproclaim

\demo{Proof} As the proof of Lemma 2.3.2 in \cite{4}.
\enddemo

Let us finish the proof of Theorem 2.2.2. We have the Mukai identification 
(it is defined by \thetag{2.1.5}) of
the transcendental periods
$$
(T(X),\,H^{2,0}(X))=(T(Y),\,H^{2,0}(Y)).
\tag{2.3.15}
$$
For general $X$ with the Picard lattice $N(X)$, it is the unique
isomorphism of the transcendental periods up to multiplication by $\pm 1$.
If $X\cong Y$, this (up to $\pm 1$) isomorphism can be extended to
$\phi:H^2(X, \bz)\cong H^2(Y, \bz)$. The restriction of $\phi$ on $N(X)$
gives then isomorphism $\phi_1:N(X)\cong N(Y)$ which is $\pm 1$ on
$A_{N(X)}=A_{N(Y)}$ under the identification \thetag{2.3.12}.
The element $\widetilde{h}=(\phi_1)^{-1}(h)$ and $f=\phi^{-1}$
satisfy Theorem 2.2.2 by Lemma 2.3.2.

The other way round, under conditions of Theorem 2.2.2, by Lemma 2.3.2,
one can construct an isomorphism $\phi_1:N(X)\cong N(Y)$ which is
$\pm 1$ on $A_{N(X)}=A_{N(Y)}$. It can be extended to
be $\pm 1$ on the transcendental periods under the Mukai identification
\thetag{2.3.15}. Then it is defined by the isomorphism
$\phi:H^2(X,\,\bz)\to H^2(Y,\,\bz)$. Multiplying $\phi$ by $\pm 1$ and
by reflections from $W^{(-2)}(N(X))$, if necessary (the group
$W^{(-2)}(N(X))$ acts identically on the discriminant group
$N(X)^\ast/N(X)$), we can assume that $\phi$ maps the
K\"ahler cone of $X$ to the K\"ahler cone of $Y$.
By global Torelli Theorem for K3 surfaces \cite{12}, 
it is then defined by an isomorphism of $X$ and $Y$.
\enddemo

\head
3. The case of Picard number 2
\endhead

\subhead
3.1. Description of general $(X,H)$ with $\rho (X)=2$ such 
that $Y\cong X$ and $\gamma (H)=1$ for $a$ odd 
\endsubhead
Here we apply results of Sect. 2 to $X$ and $Y$ with
Picard number 2.

We start with some preliminary considerations on
K3 surfaces $X$ with Picard number 2 and a primitive polarization $H$ of
degree $H^2=2a^2$, $a \ge 1$.  Thus, we assume that $\rk N(X)=2$. 
Additionally, we assume that $\gamma (H)=1$ for $H\in N(X)$
(we have this condition, if $a$ is even and $Y\cong X$). 
Let
$$
K(H)=H^\perp_{N(X)}=\bz \delta
$$
and $\delta^2=-t$ where $t>0$ is even. The $\delta\in N(X)$ is defined 
uniquely up to $\pm \delta$. 
It then follows that
$$
N(X)=[\bz H,\, \bz \delta,\, \mu H^\ast+{\delta \over 2a^2}]
$$
where $H^\ast=H/2a^2$ and $g.c.d(\mu,2a^2)=1$. 
The element    
$$
\pm \mu \mod 2a^2 \in   (\bz/2a^2)^\ast  
$$
is {\it the invariant of the pair $(N(X),H)$} 
up to isomorphisms of lattices 
with the primitive vector $H$ of $H^2=2a^2$. If $\delta $ changes to 
$-\delta$, the $\mu\mod 2a^2$ changes to $-\mu\mod 2a^2$.   
We have 
$$
(\mu H^\ast+{\delta\over 2a^2})^2=
\frac{1}{2a^2}(\mu^2-\frac{t}{2a^2}) \equiv 0\ \mod 2.
$$
Then $t=2a^2 d$, for some $d \in \bn$ and 
$\mu^2 \equiv\ d \mod 4a^2$. Thus, $d\mod 4a^2 \in (\bz/4a^2)^{\ast\,2}$. 
Obviously, $-d=\det (N(X))$.  

Any element $z\in N(X)$ can be written as 
$z=(xH+y\delta)/2a^2\ \text{where}\ x \equiv \mu y \mod 2a^2$.  
In these considerations, one can replace $H$ by
any primitive element of $N(X)$ with square $2a^2$. Thus, we have:

\proclaim{Proposition 3.1.1} Let $X$ be a K3 surface with the Picard
number $\rho = 2$ and a primitive
polarization $H$ of degree $H^2=2a^2$, $a>0$, and 
$\gamma (H)=1$ for $H\in N(X)$. 

The pair $(N(X),H)$ has the invariants $d\in \bn$ and 
$\pm \mu\mod 2a^2\in (\bz/2a^2)^\ast$ 
such that $\mu^2\equiv d \mod 4a^2$. (It follows, that $d\equiv 1\mod 4$.)

For the invariants $d$, $\mu$ we have: $\det N(X)=-d$, and  
$K(H)=H^\perp_{N(X)}=\bz \delta$ where $\delta^2=-2a^2d$. Moreover,  
$$
N(X)=[H,\delta,(\mu H+\delta)/2a^2], 
\tag{3.1.1}
$$ 
$$
N(X)=\{z=(xH+y\delta)/2a^2 \ |\ x,y\in \bz\ \text{and}\ 
x \equiv \mu y \mod 2a^2 \},
\tag{3.1.2}
$$
and $z^2=(x^2-dy^2)/2a^2$. 

For any primitive element $H^\prime \in N(X)$ with
$(H^\prime)^2 =H^2=2a^2$ and the same invariant $\pm \mu$, 
there exists an automorphism
$\phi\in O(N(X))$ such that $\phi(H)=H^\prime$.
\endproclaim

Applying Proposition 3.1.1 to $a=1$, we get that the 
pair $(N(Y), h)$ with $h^2=2$ and $\gamma (h)=1$ 
is defined by its determinant $\det N(Y)=-d$ where $d\equiv 1\mod 4$. 
Thus, from Propositions 3.1.1, we get 

\proclaim{Proposition 3.1.2} Under conditions and notations of
Propositions 3.1.1, all elements
$h^\prime=(xH+y\delta)/(2a^2) \in N(X)$ with
$(h^\prime)^2=2$ are in one to one
correspondence with integral solutions $(x,y)$ of the equation
$$
x^2-dy^2=4a^2
\tag{3.1.3}
$$
such that $x \equiv \mu y \mod 2a^2$.

The Picard lattices of $X$ and $Y$ are isomorphic, $N(X)\cong N(Y)$,
if and only if there exists such a solution. 
\endproclaim

\demo{Proof} It follows from the fact that $\gamma (h)=1$ for 
$h\in N(Y)$ since \thetag{2.3.9}. 
\enddemo   

The crucial statement is

\proclaim{Theorem 3.1.3}
Let $X$ be a K3 surface, $\rho (X)=2$ and $H$ a primitive polarization of 
$X$ of degree $H^2=2a^2$, $a>1$. 
Let $Y$ be the moduli space
of sheaves on $X$ with the Mukai vector $v=(a,H,a)$ and the canonical
$nef$ element $h=(-1,0,1)\mod \bz v$. 
Assume that $\gamma (H)=1$ (for even $a$ 
only for this case we may have $Y\cong X$). 
Then we can introduce the invariants 
$\pm \mu \mod 2a^2 \in (\bz/2a^2)^\ast$ and $d\in \bn$ of $(N(X), H)$ 
as in Proposition 3.1.1. 
Thus, we have  
$$
\gamma(H)=1,\ \det N(X)=-d\ \text{where\ } \ \mu^2 \equiv d\mod 4a^2. 
\tag{3.1.4}
$$
With notations of Propositions 3.1.1, 
all elements $\widetilde{h}=(xH+y\delta)/(2a^2) \in N(X)$
with square
$\widetilde{h}^2=2$ satisfying Theorem 2.2.2 are in one to one 
correspondence with integral solutions $(x,y)$ of the equation 
$$
x^2-dy^2=4a^2 
\tag{3.1.5}
$$
with $x \equiv\mu y \mod 2a^2$ and $x\equiv \pm 2a \mod d$. 

In particular (by Theorem 2.2.2), 
for a general $X$ with $\rho(X)=2$ and $\gamma (H)=1$ for odd $a$, 
we have $Y\cong X$ if and only if the equation $x^2-dy^2=4a^2$ 
has an integral solution $(x,y)$ with $x\equiv \mu y\mod 2a^2$ and 
$x\equiv \pm 2a \mod d$.
Moreover, a nef element $h=(xH+y\delta)/2a^2$ with $h^2=2$ defines 
the structure of a double plane on $X$ which is isomorphic to the
double plane $Y$ if and only if $x\equiv \pm 2a\mod d$. 
\endproclaim

\demo{Proof} 
Let $\widetilde{h}\in N(X)$ satisfies conditions of Theorem 2.2.2.

By Proposition 3.1.2, all primitive
$$
\widetilde{h}=(xH+y\delta)/(2a^2) \in N(X)
\tag{3.1.6}
$$
with $(\widetilde{h})^2=2$ are in one to one
correspondence with integral $(x,y)$ which satisfy the equation
$x^2-dy^2=4a^2$ and $x \equiv \mu y \mod 2a^2$,
and any integral solution of the equation $x^2-dy^2=4a^2$ with 
$x \equiv \mu y \mod 4a^2$ gives such $\widetilde{h}$.  

Let $k=aH+b\delta \in \widetilde{h}^\perp =\bz \alpha$. Then
$(k,h)=ax-byd=0$ and $(a,b)=\lambda(yd,x)$. Hence, we have
$$
(\lambda(ydH+x\delta))^2=\lambda^2(2a^2y^2d^2-2a^2dx^2)=
2a^2\lambda^2d(y^2d-x^2)=-2(2a^2)^2\lambda^2d. 
$$
Since $\alpha^2=-2d$, we get $\lambda=1/2a^2$ and
$\alpha=(ydH+x\delta)/2a^2$. There exists a unique (up to $\pm 1$)
embedding $f:K(H)=\bz\delta\to K(h)=\bz\alpha$
of one-dimensional lattices. It is given by $f(\delta)=a\alpha$
up to $\pm 1$. Thus, its dual is defined by
$f^\ast(\alpha^\ast)=a\delta^\ast$ where $\alpha^\ast=\alpha/2d$ and
$\delta^\ast=\delta/2a^2d$. To satisfy
conditions of Theorem 2.2.2, we should have
$$
\kappa(\widetilde{h})(\alpha^\ast)=\pm a \kappa(H)(\delta^\ast).
\tag{3.1.7}
$$
Further we denote $\nu=\mu^{-1}\mod 2a^2$. 
We have $u^\ast=\nu d\delta^\ast$, 
$w^\ast=a f(u^\ast)=\nu \frac{\alpha}2$,
and
$$
\kappa(\widetilde{h})(\alpha^\ast)=(-2\alpha^\ast\cdot w^\ast)
\widetilde{h}^\ast+
\alpha^\ast+N(X)
\tag{3.1.8}
$$
by \thetag{2.2.8}. Here $\widetilde{h}^\ast=\widetilde{h}/2$.
We then have $\alpha^\ast\cdot w^\ast=-\frac{\nu}2$, and
$$
\kappa(\widetilde{h})(\alpha^\ast)=
\nu \widetilde{h}^\ast+\alpha^\ast+N(X)=
(\frac{\nu x+y}2) H^\ast+(\frac{\nu yd+x}2) \delta^\ast
\tag{3.1.9}
$$
where $H^\ast=H/2a^2$. \par
\noindent We have
$u^\ast=\nu d\delta^\ast=\nu \delta/2a^2$.
By \thetag{2.2.8},
$\kappa(H)(\delta^\ast)=
(-2a^2 \delta^\ast\cdot u^\ast)H^\ast+\delta^\ast+N(X)$.
We have $\delta^\ast\cdot u^\ast=-\nu/2a^2$. It follows,
$$
\kappa (H)(\delta^\ast)=\nu H^\ast+\delta^\ast.
\tag{3.1.10}
$$
By \thetag{3.1.9} and \thetag{3.1.10}, we then get that
$\kappa (H)(a\delta^\ast)=\pm \kappa (\widetilde{h})(\alpha^\ast)$ is
equivalent to $(\nu yd+x)/2 \equiv \pm a \mod d$ and hence
$x+\nu y d \equiv \ \pm 2a \mod d$ since
the group $N(X)^\ast/N(X)$ is cyclic of order $d$ and it is generated by
$\nu H^\ast+\delta^\ast+N(X)$. Thus, finally we get
$x\equiv \pm 2a \mod d$.

This finishes the proof.
\enddemo

By Theorems 3.1.3, for given $a>1$, 
$\pm \mu \in  (\bz/2a^2)^\ast$ and $d\in \bn$ such 
that $d\equiv \mu^2 \mod 4a^2$, we
should look for all integral $(x,y)$ such that.  
$$
x^2-dy^2=4a^2,\ x\equiv \mu y \mod 2a^2,\ x\equiv \pm 2a\mod d. 
\tag{3.1.11}
$$

We first describe the set of integral $(x,y)$ such that 
$$
x^2-dy^2=4a^2,\ x\equiv \pm 2a\mod d.  
\tag{3.1.12}
$$
Considering $\pm (x,y)$, we can assume that  
$x\equiv 2a \mod d$. Then $x=2a-kd$ where $k\in \bz$. We have 
$4a^2-4akd+k^2d^2-dy^2=4a^2$. Thus, 
$$
d={y^2+4ak\over k^2}.
$$
Let $l$ be prime. Like in \cite{4}, it is easy to see that if 
$l^{2t+1}\vert k$ and $l^{2t+2}$ does not divide $k$, 
then $l^{2t+2}\vert 4ak$. 
It follows that $k=-\alpha q^2$ where $q\in \bz$, $\alpha\in \bz$ 
is square-free and 
$\alpha \vert 2a$. Then 
$$
d={y^2-4a\alpha q^2\over \alpha^2q^4}.
$$ 
It follows, 
$\alpha q \vert y$ and $y=\alpha qp$ where $p\in \bz$. We then get 
$$
d={p^2-4a/\alpha\over q^2}. 
\tag{3.1.13}
$$ 
Equivalently,  
$$
\alpha | 2a\ \text{is square-free,}\ \ p^2-dq^2={4a\over \alpha}\ . 
\tag{3.1.14}
$$
We remind that $\alpha$ 
can be negative. 

Thus, solutions $\alpha$, $(p,q)$ of \thetag{3.1.14} give all solutions  
$$
(x,y)=\pm (2a+\alpha dq^2,\ \alpha pq) 
\tag{3.1.15}
$$
of \thetag{3.1.12}. We call them as {\it associated solutions.} Thus, 
all solutions $(x,y)$ of \thetag{3.1.12} are associated solutions 
\thetag{3.1.15} to all solutions $\alpha$, $(p,q)$ of \thetag{3.1.14}.  
If one additionally assumes that $q>0$, then $(x,y)$ and 
$\alpha$, $(p,q)$ are in one to one correspondence (by our 
construction). 

Now let us consider associated solutions \thetag{3.1.15} which satisfy 
the additional condition $x\equiv \mu y\mod 2a^2$. We have 
$$
2a+\alpha dq^2\equiv \mu\alpha pq\mod 2a^2.
\tag{3.1.16}
$$
Using 
$d\equiv \mu^2\mod 4a^2$, we get 
$$
{2a\over \alpha}\equiv \mu q(p-\mu q)\mod 2a^2/\alpha.
\tag{3.1.17}
$$
Since $p^2-dq^2=4a/\alpha$, we get from {3.1.16} 
$$
-{2a\over \alpha}\equiv -p(p-\mu q)\mod 2a^2/\alpha .
\tag{3.1.18}
$$
Taking sum, we get 
$$
(p-\mu q)^2\equiv 0\mod 2a^2/\alpha.
\tag{3.1.19}
$$
Since $\alpha |2a^2$ and $\alpha$ is square-free, it follows easily that 
$a | (p-\mu q)$ and $(p-\mu q)/a$ is an integer. From 
\thetag{3.1.18}, we then get 
$$
a\vert p-\mu q\ \ \text{and}\ \ 
2\equiv \alpha p \left({p-\mu q\over a}\right)\mod 2a.
\tag{3.1.20}
$$
The condition $x\equiv \mu y\mod 2a^2$ is equivalent to \thetag{3.1.20}. 
 
From \thetag{3.1.20}, we get $\alpha | 2$. Thus,  
$\alpha =\pm 1$ or $\alpha =\pm 2$. Let us consider both cases. 

Assume $\alpha =\pm 1$. By \thetag{3.1.19}, then $2a|(p-\mu q)$, and 
we can rewrite \thetag{3.1.20} as 
$$
2a\vert p-\mu q\ \ \text{and}\ \ 
\pm 1\equiv p \left({p-\mu q\over 2a}\right)\mod a.
\tag{3.1.21}
$$
For $\alpha=\pm 1$ we have $p^2-dq^2=\pm 4a$. It follows 
$(p-\mu q)(p+\mu q)\equiv \pm 4a \mod 4a^2$. If $2a \vert (p-\mu q)$, 
then $((p-\mu q)/2a)(p+\mu q)\equiv ((p-\mu q)/2a)2p\equiv \pm 2\mod 2a$  
which is equivalent to \thetag{3.1.21}. 

Thus, for $\alpha=\pm 1$, associated solutions $(x,y)$ to 
$\alpha=\pm 1$, $(p,q)$ satisfy the additional condition 
$x\equiv \mu y\mod 2a^2$ if and only if $p\equiv \mu q\mod 2 a$. 
Equivalently, $ap\equiv \mu aq \mod 2 a^2$ which is equivalent to 
$h_1=(apH+aq\delta)/2a^2 \in N(X)$. The equation $p^2-dq^2=\pm 4a$ 
is equivalent to $(h_1)^2=\pm 2a$. We also have 
$h_1\cdot H=ap\equiv 0\mod a$. Vice versa, assume 
$h_1=(uH+v\delta)/2a^2 \in N(X)$, and $(h_1)^2=\pm 2a$, 
$h_1\cdot H\equiv 0\mod a$. We then have $h_1\cdot H=u\equiv 0\mod a$. 
Thus, $u=ap$ and $h_1=(ap+v\delta)/2a^2$. Since $h_1\in N(X)$, 
then $ap\equiv \mu v\mod 2a^2$. Thus, $v=aq$ and $p\equiv \mu q\mod 2a$. 
Since $(h_1)^2=\pm 2a$, we get $p^2-dq^2=\pm 4a$. We also remark that 
from the conditions \thetag{3.1.20} and $p^2-dq^2=\pm 4a$, 
it follows that  
$$
\text{g.c.d}(\alpha a,p)=\text{g.c.d}(\alpha a,q)=1;\ 
\text{g.c.d}(p,q)\vert (2/\alpha).
\tag{3.1.22}
$$ 
Thus, $(p,q)$ is an ``almost primitive'' solution of the equation 
$p^2-dq^2=\pm 4a$. It is primitive, if $a$ is even, and 
$\text{g.c.d}(p,q)\vert 2$, if $a$ is odd.

Now assume that $\alpha = \pm 2$. Then \thetag{3.1.20} is equivalent to 
$$
a\vert p-\mu q\ \ \text{and}\ \ 
\pm 1\equiv p \left({p-\mu q\over a}\right)\mod a.
\tag{3.1.23}
$$
Assume that $a\vert p-\mu q$ and $p^2-dq^2=\pm 2 a$. 
Then  $(p-\mu q)(p+\mu q)\equiv \pm 2a\mod 4a^2$ and 
$((p-\mu q)/a)(p+\mu q)\equiv ((p-\mu q)/a)2p\equiv \pm 2 \mod a$. 
If $a$ is odd, this is equivalent to \thetag{3.1.23}. 

Assume that $a$ is even. If $(p-\mu q)/a$ is even, then $p+\mu q$ is 
also even. From $p^2-dq^2=\pm 2 a$, we then get  
$((p-\mu q)/2a)(p+\mu q)\equiv \pm 1\mod 2a$ which is impossible for 
even $p+\mu q$. 
Thus, $p-\mu q \equiv a \mod 2a$ and $\mu q \equiv p+a\mod 2a$. 
From $((p-\mu q)/a)(p+\mu q)\equiv \pm 2 \mod 4a$, we then get 
$((p-\mu q)/a)(2p+a)\equiv \pm 2 \mod 2a$ and 
$((p-\mu q)/a)p + ((p-\mu q)/a)(a/2) \equiv \pm 1 \mod a$. Since 
$(p-\mu q)/a$ is odd, it follows that \thetag{3.1.23} never satisfies  
for even $a$.

Thus, we get that for $\alpha = \pm 2$ the number $a$ is 
odd and the condition $x\equiv \mu y\mod 2a^2$ is equivalent to 
$p\equiv \mu q\mod a$. Let us consider this case. 
From $p^2-dq^2=\pm 2a$ and $d$ odd, we get 
that $p\equiv q\mod 2$. Since $a$ is odd, we then get 
$p\equiv \mu q\mod 2a$ and $p^2\equiv \mu^2q^2\mod 4a$. 
This contradicts $p^2-dq^2=\pm 2a$ because $d\equiv \mu^2\mod 4a^2$. 
Thus, $\alpha=\pm 2$ is impossible for odd $a$ too.  

Finally we get the main results. 

\proclaim{Theorem 3.1.4} 
With conditions of Theorem 3.1.3, 
for a general $X$ with $\rho (X)=2$ and $\gamma (H)=1$ for odd a, 
we have $Y\cong X$, if and only if at least for one $\alpha = \pm 1$ 
there exists integral $(p,q)$ such that 
$$
p^2-dq^2={4a\over \alpha}\ \ \text{and}\ \ 
p\equiv \mu q\mod {2a }.
\tag{3.1.24}
$$
Solutions $(p,q)$ of \thetag{3.1.24} 
are ``almost primitive'', they satisfy \thetag{3.1.22}. 

Solutions $(p,q)$ of \thetag{3.1.24} give all solutions \thetag{3.1.5} 
of Theorem 3.1.3 as associated solutions 
$$
(x,y)=\pm (2a+\alpha d q^2,\alpha pq).
$$
\endproclaim

Interpreting, like above, solutions $(p,q)$ of \thetag{3.1.24} as elements of 
$N(X)$, we also get 

\proclaim{Theorem 3.1.5} 
With conditions of Theorem 3.1.3, for a general $X$ with 
$\rho (X)=2$ and $\gamma (H)=1$ for odd $a$,  
we have $Y\cong X$, if and only if at least for one $\alpha= \pm 1$ 
there exists $h_1\in N(X)$ such that 
$$
h_1^2=2\alpha a\ \ \text{and}\ \ h_1\cdot H\equiv 0 \mod a.
$$
\endproclaim

Applying additionally Theorem 2.2.2, we get the following simple 
sufficient condition when $Y\cong X$ which is valid for $X$ with 
any $\rho (X)$. This is one of the main results of the paper. 

\proclaim{Theorem 3.1.6}
Let $X$ be a K3 surface and $H$ a primitive polarization of degree 
$2a^2$, $a \geq 2$. Let $Y$ be
the moduli space of sheaves on $X$ with the Mukai vector $v=(a,H,a)$.

Then $Y\cong X$ if at least for one $\alpha =\pm 1$   
there exists $h_1\in N(X)$ such that 
$$
(h_1)^2=2\alpha a,\ \  h_1\cdot H\equiv 0\mod a ,
\tag{3.1.25}
$$
and $\gamma (H)=1$ for $H\in [H,h_1]_{\text{\pr}}$ where 
$[H,h_1]_{\text{\pr}}$ is the primitive sublattice of $N(X)$  
generated by $H,h_1$.

This condition is necessary to have $Y\cong X$ and $\gamma (H)=1$ for 
$a$ odd, if either $\rho(X)=1$, or $\rho(X)=2$,  
and $X$ is a general K3 surface with its Picard lattice
(i. e. the automorphism group of the transcendental periods
$(T(X),H^{2,0}(X))$ is $\pm 1$).
\endproclaim

\demo{Proof} The cases $\rho (X)\le 2$ had been considered.
We can assume that $\rho (X)>2$. Let $N=[H,h_1]_{\pr}$. All considerations
above for $N(X)$ of $\rk N(X)=2$ will be valid for $N$.
We can construct an associated with $h_1$ solution
$\widetilde{h}\in N$ with $\widetilde{h}^2=2$ such that $H$ and
$\widetilde{h}$ satisfy conditions of Theorem 2.2.2 for $N(X)$
replaced by $N$. It is easy to see that the conditions (b) and (c)
will be still satisfied if we extend $f$ in (b) $\pm$ identically on
the orthogonal complement $N^\perp_{N(X)}$. 
\enddemo 

It seems, many known examples of $Y\cong X$ (e. g. see 
\cite{2}, \cite{8}, \cite{17}) follow from Theorem 3.1.6. 
Theorems 3.1.5 and 3.1.6 are also interesting because they give 
a very clear geometric interpretation of some elements 
$h_1\in N(X)$ with negative square $(h_1)^2$ (for negative $\alpha$).

Below we consider an application of Theorem 3.1.4.    

\subhead
3.2. Description of all divisorial conditions on moduli $(X,H)$ 
such that $Y\cong X$, and $\gamma (H)=1$ for odd $a$ 
\endsubhead
Further we use the following notations. We fix $a\in \bn$, 
$\alpha \in \{1,-1\}$    
and $\overline{\mu}=\{\mu,-\mu\}\subset (\bz/ 2a^2)^\ast$. 
We denote by $\Da(a)^{\overline{\mu}}_\alpha$ the set of all $d\in \bn$ 
such that $d\equiv \mu^2 \mod 4a^2$ and there exists an integral 
$(p,q)$ such that $p^2-dq^2=4a/\alpha$ and 
$p\equiv \mu q\mod 2a $. We denote by 
$\Da(a)^{\overline{\mu}}$ the union of $\Da(a)^{\overline{\mu}}_\alpha$ 
for all $\alpha \in \{1,-1\}$, by   
$\Da(a)_\alpha$ the union of $\Da(a)^{\overline{\mu}}_\alpha$ 
for all $\overline{\mu}=\{\mu,-\mu\}\subset (\bz/ 2a^2)^\ast$, 
and by $\Da(a)$ the union of all $\Da(a)_\alpha$ for all 
$\alpha \in \{1,-1\}$. 

Assume that $X$ is a K3 surface with a primitive polarization $H$ 
of degree $H^2=2a^2$ where $a>1$. The moduli space $Y$ of sheaves 
on $X$ with Mukai vector $v=(a,H,a)$ has the canonical $nef$ 
element $h=(-1,0,1) \mod \bz v$ with $h^2=2$. It follows that 
$Y$ is never isomorphic to $X$ if $\rho (X)=1$. Since the dimension of 
moduli of $(X,H)$ is equal to $20-\rho (X)$, it follows that describing 
general $(X,H)$ with $\rho (X)=2$ and $Y\cong X$, we at the same time 
describe all possible divisorial conditions on moduli of $(X,H)$ 
when $Y\cong X$. See \cite{9}, \cite{10} and also \cite{3}. 
They are described by invariants of the pairs $(N(X),H)$ 
where $\rk N(X)=2$. By Theorem 3.1.4, we get    

\proclaim{Theorem 3.2.1} All possible divisorial conditions on 
moduli of polarized K3 surfaces $(X,H)$ with a primitive polarization 
$H$ with $H^2=2a^2$, $a>1$, which imply $Y\cong X$ and $\gamma (H)=1$ for 
$a$ odd, are labelled by 
the set ${\Cal Div}(a)$ of all pairs  
$$
(d,\,\overline{\mu})
$$
where $\overline{\mu}=\{\mu,\, -\mu\}\subset (\bz/2a^2)^\ast$,  
$d\in \Da(a)^{\overline{\mu}}=
\bigcup_{\alpha}{\Da(a)^{\overline{\mu}}_\alpha}$. Here 
$\alpha \in \{1,-1\}$.  

For any $\overline{\mu}=\{\mu,\, -\mu\}\subset (\bz/2a^2)^\ast$  
and any $\alpha \in \{1,-1\}$ 
the set 
$$
\split
\Da(a)^{\overline{\mu}}_\alpha &=\{d={p^2-4a/\alpha\over q^2}\in \bn\ |\ 
q\in \bn,\ p\equiv \mu q\mod 2a ,\   
d\equiv \mu^2\mod 4a^2 \}\\
&\supset  
\{\left(\mu +t(2a/\alpha)\right)^2-4a/\alpha\in \bn \ |\ t\mu \equiv 1\mod a\}
\endsplit
$$
is infinite (put $q=1$, to get the last infinite subset). 

In particular, for any $a>1$, the set of possible divisorial conditions 
on moduli of $(X,H)$ which imply $Y\cong X$, is infinite. 
\endproclaim

To enumerate the sets $\Da(a)^{\overline{\mu}}_\alpha$, 
$\alpha \in \{1,-1\}$, it is 
the most important to enumerate the sets $\Da(a)_\alpha$. This 
is almost equivalent to finding all possible $d\in \bn$ such that 
\newline 
$d\mod 4a^2\in (\bz/4a^2)^{\ast 2}$ and the equation 
$p^2-dq^2=4a/\alpha$ has a solution $(p,q)$ satisfying \thetag{3.1.22} 
(it is ``almost primitive''). Each such a solution 
$(p,q)$ defines a unique (if it exists) 
$\mu\mod 2a^2$ such that $\mu^2\equiv d\mod 4a^2$ and 
$p\equiv \mu q\mod 2a $.  
The pair $(d,\,\overline{\mu})$ gives then an element of the set 
${\Cal Div}(a)$. Thus, to find all possible $\overline{\mu}$ (for 
the given $\alpha$, $d$), it is enough to find 
all  ``almost primitive'' solutions 
$(p,q)$ of the equation $p^2-dq^2=4a/\alpha$ and $\nu \mod 2a $ 
such that $p\equiv \nu q\mod 2a $.  
 
Really, let $(p,q)$ be a   
solution of the equation $p^2-dq^2=4a/\alpha$. For example, assume 
that it is primitive.   
Let $\nu\equiv p/q\mod 2a$. Then $\nu^2\equiv 
d\mod 4a $. We have $\mu=\nu+k2a \mod 2a^2$, $k\in \bz$, and 
$\mu^2\equiv \nu^2+4k\nu a   + 4k^2 a^2  
\equiv d \mod 4a^2$. 
Equivalently, 
$$
4\nu k a \equiv d-\nu^2\mod 4a^2.
$$
Finally, we get 
$$
\nu k \equiv {d-\nu^2\over 4a }\mod a
$$
which determines $\mu \mod 2a^2$ uniquely. 
If $\text{g.c.d}(p,q)=2$, then $a$ is odd. 
For this case, one should again start with $\nu \mod 2a$ such 
that $p\equiv \nu q\mod 2a $ and $\nu^2\equiv d\mod 4a$. Then 
there exists a unique lifting $\mu\mod 2a^2$ such that 
$\mu\equiv \nu\mod 2a$ and $\mu^2\equiv d\mod 4a^2$.

\enddocument

\end